\documentclass[12pt]{amsart}
\usepackage{amsfonts}
\usepackage{amssymb}
\theoremstyle{definition}

\theoremstyle{remark}

\newcommand{\C}{\mathbb C}
\newcommand{\D}{\mathbb D}
\newcommand{\Sph}{\mathbb S}

\newcommand{\ds}{\displaystyle}

\begin{document}

\centerline{\bf Comptes rendus de l'Academie bulgare des Sciences}

\centerline{\it Tome 35, No 3, 1982}

\vspace{0.6in}
\begin{flushright}
{\it MATH\'EMATIQUES
\\ G\'eometrie diff\'erentielle}
\end{flushright}

\vspace{0.3in}

\centerline{\large\bf SCHUR'S THEOEM OF ANTIHOLOMORPHIC TYPE }

\vspace{0.1cm}
\centerline{\large\bf FOR QUASI-K\"AHLER MANIFOLDS}

\vspace{0.2in}
\centerline{\bf G. T. Ganchev, O. T. Kassabov}

\vspace{0.1in}
\centerline{\it (Submited by Academician B. Petkanchin on November 25, 1981)}

\vspace{0.1in}

Let $M$ be an almost Hermitian manifold with dim\,$M=2n$, a metric tensor $g$ and
an almost complex structure $J$. The Ricci tensor 
and the scalar curvature with regard to the curvature tensor $R$ are denoted 
by $S$ and $\tau$ respectively. If $\alpha$ is a 2-plane in the tangential 
space $T_pM$ with an orthonormal basis $\{x,y\}$, then its sectional curvature
$K(\alpha,p)$ is given by $K(\alpha,p)=R(x,y,y,x)$. A 2-plane $\alpha$ in 
$T_pM$ is said to be holomorphic (antiholomorphic) if $J\alpha=\alpha$ 
($J\alpha\perp\alpha$). The manifold $M$ is said to be  of pointwise
constant holomorphic (antiholomorphic) sectional curvature if the  
sectional curvature $K(\alpha,p)$ of an arbitrary holomorphic (antiholomorphic) 
2-plane $\alpha$ in $T_pM$  does not depend on $\alpha$
for every $p\in M$. The classes of K\"ahlerian $(K-)$, nearly-K\"ahlerian ($NK-$)
and quasi-K\"ahlerian ($QK-$) manifolds are characterized by the following 
identities for the covariant derivative of $J$: $(\nabla_XJ)Y=0$, $(\nabla_XJ)X=0$,
$(\nabla_{JX}J)Y=-J(\nabla_XJ)Y$ respectively. Every $K$-manifold is an $NK$-manifold 
and every $NK$-manifold is a $QK$-manifold. The curvature tensor of an $NK$-manifold 
satisfies the identity  \cite{G2}
$$
	R(x,y,z,u)=R(x,y,Jz,Ju)+R(x,Jy,z,Ju)+R(Jx,y,z,Ju) \ .  \leqno (1)
$$
This equality implies
$$
	R(x,y,z,u)=R(Jx,Jy,Jz,Ju)   \ .   \leqno (2)
$$
A $QK$-manifold whose curvature tensor satisfies (1) is said to be a $QK_2$-manifold. 
An almost Hermitian manifold whose curvature tensor satisfies (2) is said to be an $AK_3$-manifold.
For a $K$-manifold the following analogues of the classical theorem of F. Schur hold \cite{KN, ChO}: 

Let $M$ be a connected $2n-$dimensional $K$-manifold of pointwise constant holomorphic
(antiholomorphic) sectional curvature $c$. If $n\ge 2$ ($n\ge 3$), then $c$ is a global constant.
Moreover $M$ is of constant antiholomorphic (holomorphic) sectional curvature $c/4\, (c)$. 

Naveira and Hervella proved \cite{NH} that if a connected $NK$-manifold of dimension
$2n\ge 4$ has a pointwise constant holomorphic sectional curvature $\mu$, then $\mu$
is a constant and Gray proved in \cite{G1} that such a manifold is locally isometric to
one of the following manifolds: $\C^n$, $\C\mathbb P^n$, $\C\D^n$, $\Sph^6$. We proved 
in \cite{GK} that if $M$ is a connected $NK$-manifold of pointwise constant antiholomorphic
sectional curvature $\nu$ and dim\,$M\ge6$, then $\nu$ is a constant and $M$ is locally
isometric to $\C^n$, $\C\mathbb P^n$, $\C\D^n$ or $\Sph^6$.

Gray and Vanhecke proved in \cite{GV} that if a connected $QK_2$-manifold has a
pointwise constant holomorphic sectional curvature $\mu$, then $\mu$ is a constant.

In the present paper we consider $QK_2$-manifolds of pointwise constant antiholomorphic
sectional curvature.

Let $\{ e_1,\hdots,e_{2n} \}$ be an orthonormal basis of $T_pM$, $p\in M$, and
$$
	S'(x,y)=\sum_{i=1}^{2n} R(x,e_i,Je_i,Jy) \ ,
$$

\vspace{-0.1in}
$$
	\tau'=\sum_{i=1}^{2n} S'(e_i,e_i) \ .
$$
The tensors $R_1$, $R_2$ and $\psi$ are defined by
$$
	R_1(x,y,z,u)=g(y,z)g(x,u)-g(x,z)g(y,u)\ ;
$$
$$
	R_2(x,y,z,u)=g(Jy,z)g(Jx,u)-g(Jx,z)g(Jy,u)-2g(Jx,y)g(Jz,u) \ ;
$$
$$
	\begin{array}{r}\vspace{0.1cm}
		\psi(x,y,z,u)=g(Jy,z)S(Jx,u)-g(Jx,z)S(Jy,u)-2g(Jx,y)S(Jz,u) \ \  \\
	             +g(Jx,u)S(Jy,z)-g(Jy,u)S(Jx,z)-2g(Jz,u)S(Jx,y) \ .
	\end{array}
$$
In \cite{GK} we proved

{\bf Proposition.} Let $M$ be a $2n$-dimensional $AH_3$-manifold of pointwise constant anti\-holomorphic
sectional curvature $\nu$. Then 
$$
	R=\frac16\psi +\nu R_1-\frac{2n-1}3\nu R_2  \ ,  \leqno (3)
$$
$$
	(n+1)S-3S'=\frac{(n+1)\tau-3\tau'}{2n}\,g \ ,   \leqno (4)
$$

$$
	\nu= \frac{(2n+1)\tau-3\tau'}{8n(n^2-1)}  \ .  \leqno (5)
$$ 

In particular from (3) it follows that $R$ satisfies (1).

From the second Bianchi's identity
$$
	(\nabla_wR)(x,y,z,u)+(\nabla_xR)(y,w,z,u)+(\nabla_yR)(w,x,z,u)=0   \leqno (6)
$$ 
the next identities follow in a straightforward way:
$$
	(\nabla_xS)(y,z)-(\nabla_yS)(x,z)=\sum_{i=1}^{2n}(\nabla_{e_i}R)(x,y,z,e_i) \ , \leqno (7)
$$ 
$$
	\sum_{i=1}^{2n}(\nabla_{e_i}S)(x,e_i)=\frac12x(\tau)   \ .  \leqno (8)
$$

The first autor proved in \cite{G} that for a $QK_2$-manifold the following
analogue of (8) holds:
$$
	\sum_{i=1}^{2n}(\nabla_{e_i}S')(x,e_i)=\frac12x(\tau')   \ .  \leqno (9)
$$

{\bf Lemma.} Let $M$ be a connected $2n$-dimensional $QK_2$-manifold of
pointwise constant antiholomorphic sectional curvature. If $n\ge 2$, the
function $(n+1)\tau-3\tau'$ is a constant.

{\bf Proof.} Using (8) and (9) from (4) we obtain $(n-1)x((n+1)\tau-3\tau')=0$ for
an arbitrary vector $x$ and hence the function $(n+1)\tau-3\tau'$ is a constant.

Now let $M$ be an $AH_3$-manifold of pointwise constant antiholomorphic sectional 
curva\-ture $\nu$. Let $\{ x,y \}$ be unit orthogonal vectors and $x \perp Jy$.
Putting in (6) $w=x$, $x=y$, $y=Jy$, $z=Jy$, $u=y$ and taking into account (3),
we obtain
$$
	\begin{array}{r}\vspace{0.1in}
		4(n-1)x(\nu)=(\nabla_xS)(y,y)+(\nabla_xS)(Jy,Jy)-(\nabla_yS)(x,y)  \\
		-(\nabla_{Jy}S)(x,Jy)-S(JBy,x)-g(JBy,x)S(y,y)+2(2n-1)\nu g(JBy,x) \ ,
	\end{array}   \leqno (10)
$$
where $By=(\nabla_yJ)y+(\nabla_{Jy}J)Jy$. Putting in (7) $y=Jx$, $z=Jx$, 
and using (3) we obtain
$$
	\begin{array}{r}\vspace{0.1in}
		\ds(\nabla_xS)(Jx,Jx)-(\nabla_{Jx}S)(x,Jx)=\frac12\big\{ \frac12x(\tau)-\sum_{i=1}^{2n} S((\nabla_{e_i}J)Jx,e_i)  \\ \vspace{0.1in}
		+g(\delta F,Jx)S(x,x)+(\nabla_xS)(Jx,Jx)+S((\nabla_xJ)x,Jx) \big\}  \\
		-2(n-1)x(\nu)-(2n-1)\nu g(\delta F,Jx) \ ,
	\end{array}   \leqno (11)
$$

\vspace{-0.1cm}
\noindent
where $\ds\delta F=\sum_{i=1}^{2n}(\nabla_{e_i}J)e_i$.

Now we shall prove the following analogue of the theorem of Schur for
$QK_2$-manifolds.

{\bf Theorem.} Let $M$ be a connected $2n$-dimenssional $QK_2$-manifold of
pointwise constant antiholomorphic sectional curvature $\nu$. If $n\ge 3$, the
functions $\nu$, $\tau$ and $\tau'$ are constants.

{\bf Proof.} By the conditions of the theorem, (10) and (11) give respectively
$$
	\begin{array}{r}\vspace{0,1in}
		4(n-1)x(\nu)=(\nabla_xS)(y,y)+(\nabla_xS)(Jy,Jy)  \\
		-(\nabla_yS)(x,y)-(\nabla_{Jy}S)(x,Jy) \ ,
	\end{array}   \leqno (12)
$$
$$
	4(n-1)x(\nu)=\frac12 x(\tau)-(\nabla_xS)(Jx,Jx)+(\nabla_{Jx}S)(x,Jx) \ .  \leqno (13)
$$

Using the adapted basis $\{ u_1,\hdots,u_n;Ju_1,\hdots,Ju_n \}$, from (12) by summing
up we find
$$
	4(n-1)^2x(\nu)=\frac12 x(\tau)-(\nabla_xS)(Jx,Jx)+(\nabla_{Jx}S)(x,Jx) \ . 
$$
The last equality and (13) give $(n-1)(n-2)x(\nu)=0$ and hence $\nu$ is a constant.
Taking into account the lemma and (5), we obtain that $\tau$ and $\tau'$ are also constants.

\vspace {0.2cm}
\begin{flushright}
{\it Institute of Mathematics  \\
Bulgarian Academy of Sciences \\
Sofia, Bulgaria}
\end{flushright}

\end{document}